\documentclass[12pt]{amsart}

\numberwithin{equation}{section}

\newtheorem{theorem}[subsection]{Theorem}
\newtheorem{lemma}[subsection]{Lemma}
\newtheorem{proposition}[subsection]{Proposition}
\newtheorem{corollary}[subsection]{Corollary}
\newtheorem{conjecture}[subsection]{Conjecture}
\newtheorem{remark}[subsection]{Remark}

\theoremstyle{definition}

\newtheorem{definition}[subsection]{Definition}

\theoremstyle{remark}

\newtheorem{example}[subsection]{Example}

\newcommand{\bz}{{\bf{Z}}}

\newcommand{\zp}{{{\bf{Z}}/p}}
\newcommand{\zps}{{{\bf{Z}}/p^2}}
\newcommand{\zpr}{{{\bf{Z}}/p^r}}

\newcommand{\fp}{{\bf{F}}_p}

\newcommand{\field}{{\bf{F}}}

\newcommand{\F}[1]{{\mathcal {L}}_{#1}}
\newcommand{\B}{\mathcal {B}}

\newcommand{\fieldvg}{{\field}[V]^G}
\newcommand{\fieldv}{{\field}[V]}
\newcommand{\fvpp}{\field[V_{p+1}]}
\newcommand{\qhom}[2]{H^{#1}(Q,#2)}
\newcommand{\am}{\alpha}
\newcommand{\V}{{\mathcal {V}}}

\renewcommand{\Im}{\mathop{\rm Im}}
\newcommand{\tr}{\mathop{\rm Tr}}
\DeclareMathOperator{\Tr}{Tr}

\newcommand{\Span}{\mathop{\rm span}}

\newcommand{\lm}{\mathop{\rm LM}}

\newcommand{\del}{\Delta}

\newcommand{\Qg}{\overline{\sigma}}
\newcommand{\Ln}{N^L(x_{p+1})}
\newcommand{\binomial}[2]{\genfrac{(}{)}{0pt}{}{#1}{#2}}

\newcommand{\lato}[1]{\stackrel{#1}{\longrightarrow}}

\title[Decomposing symmetric powers]{Decomposing symmetric powers of certain modular representations of cyclic groups}

\author{R.\ James Shank}
\address{Institute of Mathematics, Statistics \&  Actuarial Science \\
 \hfil\break\indent University of Kent, Canterbury, CT2 7NF, UK}
\email{R.J.Shank@kent.ac.uk}

\author{David L.\ Wehlau}
\address{Department of Mathematics \& Computer Science \\
 \hfil\break\indent  Royal Military College, Kingston, Ontario, Canada,  K7K 7B4}
\email{wehlau@rmc.ca}

\thanks{The research of the first author is supported by grants from EPSRC.\\
\indent The research of the second author is supported by grants from ARP and NSERC}

\subjclass{13A50}
\date{\today}

\begin{document}
\begin{abstract}
For a prime number  $p$, we construct a generating set for the ring of invariants for the $p+1$ dimensional
indecomposable modular representation of a cyclic group of order $p^2$, and show that the Noether number for the representation
is $p^2+p-3$. We then use the constructed invariants
to explicitly describe the decomposition of the symmetric algebra as a module over the group ring, confirming
the Periodicity Conjecture of Ian Hughes and Gregor Kemper for this case.  In the appendix, we use our results to
compute the Hilbert series for the corresponding ring of invariants together with some other related generating
functions.
\end{abstract}

\maketitle

\begin{center}
\em
This paper is dedicated to Gerry Schwarz, on the occasion of his sixtieth
birthday.
\end{center}

\section{Introduction}

Suppose that $V$ is a finite dimensional representation of a finite group $G$ over a field $\field$,
i.e., $V$ is a finitely generated module over the group ring $\field G$.
The action of $G$ on $V$ induces an action on the dual $V^*$ which extends to an action by
algebra automorphisms on the symmetric algebra $\field[V]:=S(V^*)$.
The elements of $V^*$, and thus also the elements of $\field[V]$, represent
$\field$-valued functions on $V$.   If $\{x_1, x_2,\ldots,x_n\}$
is a basis for $V^*$ then $\field[V]$ can be identified with the ring of polynomials
$\field[x_1,x_2,\ldots,x_n]$.
Let $\field[V]_d$ denote the subspace of homogeneous polynomials
of degree $d$. Since the action of $G$ preserves degree, $\field[V]_d$ is a module over $\field G$
and $$\field[V]=\bigoplus_{d=0}^{\infty}\field[V]_d $$
is a decomposition into a direct sum of finite dimensional $\field G$-modules.
Of course $\field[V]_d$ is precisely the
the $d^{\rm th}$ symmetric power of $V^*$. Understanding the
action of $G$ on $\field[V]_d$, and hence the action on $\field[V]$, is an important problem in representation
theory. The primary goal is to write $\field[V]_d$ as a direct sum of indecomposable
$\field G$-modules, refining the given decomposition of $\field[V]$.
This means decomposing $\field[V]_d$ for infinitely many $d$.
An important aspect of the group action is the {\it ring of invariants}
$$\field[V]^G:=\{f\in\field[V]\mid g(f)=f,\ \forall g\in G\},$$
a finitely generated subalgebra of $\field[V]$. A fundamental problem in invariant theory is the
construction of a finite generating set for $\field[V]^G$. Since $G$ is finite,
$\field[V]$ is a finite module over $\field[V]^G$. Thus $\field[V]$ is a module over both
$\field[V]^G$ and $\field G$. Perhaps the right approach is to study $\field[V]$ as a finitely generated
module over the extended group ring $\field[V]^G G$. Certainly, in the work of both Karagueuzian \& Symonds \cite{ka-sy}
and Hughes \& Kemper \cite{hughes-kemper}, the finite $\field[V]^G$-module structure of $\field[V]$ has been used
to reduce decomposing $\field[V]$ over $\field G$ to a finite problem.

For the remainder of the paper, we assume that $\field$ has characteristic $p$ for a prime number $p$,
and that $G\cong \zpr$ is a cyclic group of order $p^r$. Choose a generator $\sigma$ for $G$.
The isomorphism type of a representation of $G$ is determined by the Jordan canonical form
of $\sigma$.  Since the order of $\sigma$ is a power of $p$, and since a field of characteristic $p$
has no non-trivial $p^{\rm th}$ roots of unity, all the eigenvalues of $\sigma$ must be 1.
If $m\leq p^r$, then the $m\times m$ matrix over $\field$ consisting of a single Jordan block
with eigenvalue $1$ determines an indecomposable $\field G$-module which we denote by $V_m$.
Note that if $m>p^r$, then the matrix has order greater than $p^r$ and does not determine a representation of $G$.
It follows from the form of the matrix that $V_m$ is faithful if and only if $p^{r-1}<m\leq p^r$, and that
$V_m$ is a cyclic $\field G$-module.  It is clear that if the Jordan canonical form of $\sigma$
consists of more than one Jordan block then the representation will be decomposable.
Thus the complete set of inequivalent indecomposable
 $\field G$-modules are, up to isomorphism, $V_1$, $V_2$, \ldots, $V_{p^r}$.
Furthermore, from the Jordan canonical form it is easy to see
that these modules are naturally embedded into one another:
$V_1 \subset V_2 \subset V_3 \subset \dots \subset V_{p^r}$.
Note that the one dimensional space of $G$-fixed points, $V_m^G \cong V_1$ is the {\em socle} of $V_m$.
Moreover, $V_1$ is the unique irreducible module,
$V_{p^r}\cong \field G$ is the unique projective indecomposable,
and an $\field G$-module is projective if and only if it is injective (see, for example, \cite[Ch.~II]{alperin}).
Also, it is easy to see that the representation $V_m$ is induced from a representation of a proper subgroup
of $G$ if and only if $p$ divides $m$.

For $f\in\field[V_n]$, we define the {\it norm} of $f$, denoted by $N^G(f)$, to be the product over the $G$-orbit of $f$.
Clearly $N^G(f)\in \field[V_n]^G$.
For a subgroup $L=\langle\sigma^{p^t}\rangle$, we define the {\it relative transfer} $\Tr_L^G:=\sum_{i=0}^{p^t-1}\sigma^i\in \field G$.

The two main results we prove in this article concern the representation $V_{p+1}$
and are stated as Theorem~\ref{invthm} and Theorem~\ref{decompthm0} below.
The following example illustrates these theorems.
\begin{example}\label{3-example}
  Let $\field$ be any field of characteristic $p=3$.
  We consider the indecomposable four dimensional representation $V_4$
  of the cyclic group $G=\bz/9$ of order 9.
  The group $G$ contains the subgroup $L$ of order 3.
  Theorem~\ref{invthm} asserts that $\field[V_4]^G$ is generated by
  $M = N^G(x_3) =  x_3^3  - x_3 x_2^2 + x_3^2 x_1 + x_3 x_2 x_1$,
  $N = N^G(x_4) = x_4^9  - x_4^3 x_ 3^6 + \ldots$,
  and elements from the image of the relative transfer, $\Tr_L^G(\field[x_4^3 - x_4 x_1^2, x_3,x_2,x_1])$.
  In fact, a Magma \cite{magma} computation shows that $\field[V_4]^G$ is minimally generated by
 $M$ and $N$ together with the following 9 invariants:
\begin{eqnarray*}
\Tr_L^G(x_3) &=& x_1,\\
 \Tr_L^G(-x_3^2) &=& x_2^2 + 
                              x_1 x_3 - x_1 x_2 - x_1^2,\\
 \Tr_L^G(-x_2 x_3^2) &=& x_2^3  
                                - x_1 x_2^2 - x_1^2 x_3  + x_1^3,\\
 \Tr_L^G( -x_3(x_4^3 - x_4 x_1^2)) &=& x_2 x_3^3 + \ldots,\\
 \Tr_L^G( x_3^2(x_4^3 - x_4 x_1^2)) &=& x_2 x_3^4 + \ldots,\\
 \Tr_L^G( x_2 x_3^2(x_4^3 - x_4 x_1^2)) &=& x_2^2 x_3^4 + \ldots,\\
 \Tr_L^G( -x_3(x_4^3 - x_4 x_1^2)^2) &=& x_3^7 + \ldots,\\
 \Tr_L^G( -x_3^2(x_4^3 - x_4 x_1^2)^2) &=& x_3^8 + \ldots,\\
 \Tr_L^G( x_2 x_3^2(x_4^3 - x_4 x_1^2)^2) &=& x_2 x_3^8 + \ldots\ .
 \end{eqnarray*}

   The first few homogeneous components of $\field[V_4]$ decompose into
   indecomposable $\field G$-modules as follows:
  $$
  \begin{array}{ccccrcrcr}\displaystyle
    \fieldv_0 & \cong & V_1, &    &    &      &\\   
    \fieldv_1 & \cong & V_4,    & &    &    &      &  \\  
    \fieldv_2 & \cong  & V_7 & \oplus & V_3,    &     & \\
    \fieldv_3 & \cong  & V_2 & \oplus & V_3 &\oplus  & V_6 & \oplus &V_9,\\ 
    \fieldv_4 & \cong  & V_5 & \oplus & 2\,V_3 &\oplus  & V_6 & \oplus &2\,V_9,\\ 
    \fieldv_5 & \cong  & V_8 & \oplus & 3\,V_3 &\oplus  & 2\,V_6 & \oplus &3\,V_9,\\ 
    \fieldv_6 & \cong  & &  & 4\,V_3 &\oplus  & 3\,V_6 & \oplus &6\,V_9,\\
    \fieldv_7 & \cong  & &  & 5\,V_3 &\oplus  & 4\,V_6 & \oplus &9\,V_9,\\
    \fieldv_8 & \cong  & &  & 6\,V_3 &\oplus  & 5\,V_6 & \oplus &13\,V_9 \\
    \fieldv_9 & \cong  & V_1  & \oplus &7\,V_3&\oplus &6\,V_6  &\oplus&18\,V_9\\   
    \fieldv_{10} & \cong & V_4    & \oplus &8\,V_3&\oplus &7\,V_6  &\oplus&24\,V_9\\  
    \fieldv_{11}& \cong  & V_7 & \oplus &10\,V_3&\oplus &8\,V_6  &\oplus&31\,V_9\\
    \fieldv_{12} & \cong  & V_2 & \oplus & 11\,V_3 &\oplus  & 10\,V_6 & \oplus &40\,V_9,\\ 
    \fieldv_{13} & \cong  & V_5 & \oplus & 13\,V_3 &\oplus  & 11\,V_6 & \oplus &50\,V_9,\\ 
    \fieldv_{14} & \cong  & V_8 & \oplus & 15\,V_3 &\oplus  & 13\,V_6 & \oplus &61\,V_9,\\ 
    \fieldv_{15} & \cong  & &  & 17\,V_3 &\oplus  & 15\,V_6 & \oplus &75\,V_9,\\
    \fieldv_{16} & \cong  & &  & 19\,V_3 &\oplus  & 17\,V_6 & \oplus &90\,V_9,\\
    \fieldv_{17} & \cong  & &  & 21\,V_3 &\oplus  & 19\,V_6 & \oplus &107\,V_9 .
 \end{array}
 $$
 The (one dimensional) socle of the non-induced indecomposable summand in $\field[V_4]_i$ for $i=0,1,2,\dots,5$
may be chosen to contain $1$, $x_1$, $x_1^2$, $M$, $x_1 M$ and $x_1^2 M$ respectively.

  For degrees $d \geq 9$, write $d = 9a + c$ where $0 \leq c \leq 8$.  Then
  $\fieldv_d \cong \fieldv_c \oplus \alpha\, V_3 \oplus \beta\, V_6 \oplus \gamma\, V_9$
  for some  non-negative integers $\alpha, \beta$ and $\gamma$.
  Furthermore if $0 \leq c \leq 5$ and if we denote by $f$ a non-zero element of the socle
  of the non-induced summand in $\fieldv_c$, then the non-induced summand in $\fieldv_d$
  may be chosen such that its socle is spanned by $N^a f$.
\end{example}

In Section 2 we develop tools for decomposing $\field[V_n]$ as an $\field G$-module.
We then specialise to $r=2$ and $n=p+1$.
In Section~\ref{comp_inv_sec} we construct generators for $\field[V_{p+1}]^{\zp^2}$.
We apply the ``ladder technique'' described in \cite[\S 7]{Shank-Wehlau:cmipg}, using group cohomology and
a spectral sequence argument, to
prove the following.
\begin{theorem}\label{invthm}
Suppose $G\cong\zps$ and let $L\cong\zp$ denote its non-trivial proper subgroup.
The ring of invariants $\fvpp^G$ is generated by $N^G(x_p)$, $N^G(x_{p+1})$ and
elements from the image of the relative transfer,
$\tr_L^G\left(\field[\Ln, x_p,\ldots,x_1]\right)$.
\end{theorem}

Recall that the Noether number of a representation is the largest degree of an element in a minimal homogeneous generating set
for the corresponding ring of invariants. In Section~\ref{noeth_sec} we use the generating set given by Theorem~\ref{invthm}
to show that, for $p>2$, the Noether number for $V_{p+1}$ is $p^2+p-3$.
In Section~\ref{decomp sec} we use the constructed generating set to describe the $\field\zp^2$-module structure of $\field[V_{p+1}]$,
confirming the Periodicity Conjecture of Hughes \& Kemper \cite[Conjecture~4.6]{hughes-kemper}
in this case and proving the following.
\begin{theorem}\label{decompthm0}
  Let $G\cong\zps$ and let $d$ be any non-negative integer.
    In the decomposition of $\fvpp_d$ into a direct sum of indecomposable
    $\field G$-modules there is at most one indecomposable summand $V_m$ which
     is not induced from a representation of a proper subgroup.  In particular, writing $d= a p^2 + b p + c$
     where $0 \leq b,c < p$, there is exactly one non-induced indecomposable summand
     when $b \leq p-2$ and $\fvpp_d$ is an induced module when $b=p-1$.
    Moreover,  if $b \leq p-2$ then the non-induced indecomposable summand is
    isomorphic to $V_{cp+b+1}$ and we may choose the decomposition of $\fvpp_d$ such that
    the socle of this summand, $V_{cp+b+1}^G$, is spanned by the invariant
    $N^G(x_{p+1})^a N^G(x_p)^b x_1^c$.
\end{theorem}

We note that Symonds, in a recent paper \cite{symonds} based on his joint work with Karagueuzian \cite{ka-sy}, has proven
the Periodicity Conjecture of Hughes \& Kemper.  He goes on to prove
that for $p^{r-1}<n<p^r$ and $d<p^r$, the $\field \zp^r$-module $\field[V_n]_d$ is isomorphic to $\Omega^{-d}\Lambda^d(V_{p^r-n})$
modulo induced modules \cite[Corollary~3.11]{symonds}. Here $\Lambda^d$ denotes the $d^{\rm th}$ exterior power and
$\Omega^{-d}$ denotes the  $d^{\rm th}$ cokernel of a minimal injective resolution (see \cite[page 30]{benson:rep1}).
It is instructive to compare this with Theorem~\ref{decompthm0} and Example~\ref{3-example}.

In the Appendix we compute the Hilbert series of $\fvpp^\zps$. We also compute generating functions encoding the number of
summands of each isomorphism type in $\fvpp_t$.

\section{Preliminaries}\label{prelim_sec}
  Let $G =\langle\sigma\rangle\cong \zpr$.  It will be convenient to define $\Delta:=\sigma-1\in\field G$.
   It is easy to see that $\Delta$ acts as a twisted derivation on $\field[V_n]$, i.e.,
$\Delta(a\cdot b)=a\Delta(b)+\Delta(a)\sigma(b)$.
We denote the full transfer, $\tr_{\langle1\rangle}^G$, by $\tr^G$ and
the image of the relative transfer, $\tr_L^G(\field[V_n]^L)$, by $\Im \Tr_L^G$.
Clearly $\Im \tr_L^G$ is an ideal in $\field[V_n]^G$.
A simple calculation with binomial coefficients shows that
$\Delta^{p^t}=\sigma^{p^t}-1$ and $\Delta^{p^t-1}=(\sigma^{p^t}-1)/(\sigma-1)=\Tr_L^G$.
We denote the group cohomology of $G$ with coefficients in the $\field G$-module $W$ by
$H^*(G,W)$. Note that $H^0(G,W)$ is just the fixed subspace $W^G$.
Furthermore, since $G$ is cyclic,
$H^{2i-1}(G,W)=\ker (\tr^G|_W)/\Im(\Delta|_W)$ and $H^{2i}(G,W)=\ker(\Delta|_W)/\Im (\tr^G|_W)$ for $i>0$
(see, for example, \cite[\S 2.1]{evens}).
It is clear from the definition of group cohomology that $H^i(G,P)=0$ if $i>0$ and $P$ is projective.
Thus $H^1(G,V_{p^r})=H^2(G,V_{p^r})=0$. Furthermore, if $V_m$ is generated as an $\field G$-module by
$e$, then $H^0(G,V_m)=V_m^G=\Span_{\field}(\Delta^{m-1}(e))$ and $\{e,\Delta(e),\ldots,\Delta^{m-1}(e)\}$ is a vector space
basis for $V_m$. If we identify $V_{m-1}$ with the submodule $\Delta(V_m)$, then, for $m<p^r$, $H^1(G,V_m)$
is the one dimensional vector space $V_m/V_{m-1}$ and $H^2(G,V_m)$ is the one dimensional vector space $V_m^G$.

Let  $W$ be any finite dimensional $\field G$-module.
Define $\F{t}(W) := \del^{t-1}(W)$.
Clearly $\F{i+1}(W)\subseteq\F{i}(W)$. Furthermore, since $\sigma$ has order $p^r$, $\F{p^r+1}(W)=0$.
Thus we have the following filtration of $W$ by $\field G$-modules:
$$W=\F{1}(W) \supseteq \F{2}(W) \supseteq \F{3}(W) \supseteq \dots \supseteq \F{p^r}(W)\supseteq \F{p^r+1}(W)=0.$$
This filtration of $W$ obviously induces a filtration of the subspace $W^G$:
$$W^G =\F{1}^G(W) \supseteq \F{2}^G(W) \supseteq \F{3}^G(W) \supseteq  \dots \supseteq \F{p^r}^G(W)\supseteq\F{p^r+1}^G(W)= 0$$
where $\F{t}^G(W) := \F{t}(W) \cap W^G$.

\begin{definition}
For a non-zero $f\in W$,
 we define the {\it length} of $f$, denoted by $\ell(f)$, by
$\ell(f)\geq t \iff f\in \F{t}(W)$. Note that $1\leq\ell(f)\leq p^r$.
We will refer to the above filtration of $W$ as the {\it length filtration} and say that
a basis $\B$ for $W^G$ is {\it compatible} with the length filtration if $\F{t}^G(W)\cap\B$ is
a basis for $\F{t}^G(W)$ for all $t$ (using the convention that the empty set is a basis for the
zero vector space).
\end{definition}

\begin{lemma}\label{dim lemma} If $W$ is a finite dimensional $\field G$-module, then
$$\dim(W)=\sum_{t=1}^{p^r}t\left(\dim(\F{t}^G(W))-\dim(\F{t+1}^G(W)\right).$$
\end{lemma}
\begin{proof}
Choose a decomposition of $W$ into indecomposable $\field G$-modules. For each indecomposable summand,
choose a basis in which $\sigma$ is in Jordan canonical form. The union of these bases gives a basis
for $W$. Intersecting this basis with $W^G$ gives a basis for $W^G$, say $\B$, which is compatible with the
length filtration. It is clear that the number of elements in $\B\cap(\F{t}^G(W)\setminus \F{t+1}^G(W))$
coincides with the number of indecomposable modules in the decomposition which are isomorphic to $V_t$,
giving the required formula.
\end{proof}

Suppose that $W$ is a finite dimensional $\field G$-module and $\B$ is a basis for $W^G$
which is compatible with the length filtration. For each $\alpha\in \B$ choose $\gamma\in W$
with $\del^{\ell(\alpha)-1}(\gamma)=\alpha$. (The existence of a suitable $\gamma$ follows from
the definition of length.) Define $V(\alpha)$ to be the $\field G$-module generated
by $\gamma$. Note that $\alpha$ spans the socle of $V(\alpha)$ and that $\dim(V(\alpha))=\ell(\alpha)$.

\begin{proposition}\label{decomp prop}
$$W = \bigoplus_{\alpha\in \B} V(\alpha)\ .$$
\end{proposition}
\begin{proof}
The natural homomorphism of the external direct sum of the $V(\alpha)$ to $W$ is injective
on the socle and is therefore injective. Thus the internal sum of the $V(\alpha)$ is direct.
It follows from Lemma~\ref{dim lemma} that the dimension of $W$ coincides with
the dimension of $\oplus_{\alpha\in\B}V(\alpha)$, giving equality.
\end{proof}

The above shows how we may obtain a direct sum decomposition of $W$ into indecomposable submodules
from any basis of $W^G$ which is compatible with the length filtration of $W^G$.  Clearly every
such decomposition arises in this way. Furthermore, an element $f \in W^G$ has length $t$ if
and only if there is an $\field G$ decomposition $W = W' \oplus V_t$ with $f$ spanning $V_t^G$.

  Note that if $f,h \in \fieldvg$ then $\ell(fh) \geq \ell(f)$.  To see this write $f = \del^{\ell(f)-1}(F)$.
  Then $fh = \del^{\ell(f)-1}(Fh)$.  In general it may happen that $\ell(fh) > \max\{\ell(f),\ell(h)\}$.
  Computer computations together with various results, such as Proposition~\ref{stretch prop}, lead us to make the following conjecture.
  \begin{conjecture}
    Suppose $f,h \in \fieldvg$ with $\ell(f) \equiv 0 \pmod p$.  Then $\ell(fh) \equiv 0 \pmod p$.
  \end{conjecture}


For $n\leq p^r$, choose an $\field G$-module generator $x_n$ for $V_n^*$ and define $x_i=\Delta^{n-i}(x_n)$ for $i=1,\ldots,n-1$.
Then  $\{x_1,x_2,\dots,x_n\}$ is a basis of $V_n^*$.
Let $\overline{\field}$ denote the algebraic closure of $\field$ and define
$\overline{V_n}:=\overline{\field}\otimes_{\field} V_n$.
Let $\{e_1,e_2,e_3,\dots,e_n\}$ denote the basis for $\overline{V_n}$ dual to $\{1\otimes x_1,\ldots,1\otimes x_n\}$.
Note that $e_1$ generates $\overline{V_n}$ as an $\overline{\field} G$-module and that $\Delta(e_n)=0$.
Using the inclusion $\field\subseteq \overline{\field}$, allows us to interpret elements of $\field[V_n]$
as regular functions on $\overline{V_n}$, i.e., we identify $\field[V_n]$ in a natural way with a subset
of $\overline{\field}[\overline{V_n}]$.
For a subset $X\subseteq \overline{\field}[\overline{V_n}]$, define
$\V(X)=\{v\in\overline{V_n}\mid f(v)=0\ \forall f\in X\}$.

\begin{lemma} \label{length lemma}
Suppose $p^{r-1}<n\leq p^r$ and
let $H$ denote the subgroup $\langle\sigma^{p^{t+1}}\rangle\cong\bz/{p^{r-t-1}}$ of
$G = \langle\sigma\rangle\cong \bz/{p^r}$ where $0 \leq t \leq r-1$.
\begin{enumerate}
  \item $\V(\Im \Tr_H^G) = \overline{V_n}^{\bz/p^{r-t}} =
                        \mathrm{span}_{\overline{\field}}\{e_{n-p^{t}+1},e_{n-p^{t}+2},\dots,e_{n-1},e_n\}$.
  \item For $f \in \field[V_n]^G$, if $\ell(f) \geq p^t+1$ then
$$f \in \sqrt{\Im {\rm Tr}_H^G} = \left((x_1,x_2,\dots,x_{n-p^t})\field[V_n]\right)\cap \field[V_n]^G.$$
\end{enumerate}
\end{lemma}

\begin{proof}
The equality   $V_n^{\bz/p^{r-t}}=\mathrm{span}_\field\{e_{n-p^{t}+1},e_{n-p^{t}+2},\dots,e_{n-1},e_n\}$
is easily verified. The equality  $\V(\Im \Tr_H^G) = \overline{V_n}^{\bz/p^{r-t}}$ follows from
\cite[Proposition~12.5]{fleisch} (see also \cite[Theorem 12]{chksw}). This equality of sets may be expressed
equivalently as the equality of ideals $\sqrt{\Im \Tr_H^G} = \left((x_1,x_2,\dots,x_{n-p^t})\field[V_n]\right)\cap \field[V_n]^G$
(see, for example, \cite[Proposition~11]{chksw}).
Thus it only remains to show that  if $\ell(f) \geq p^t+1$ then $f \in \sqrt{\Im \Tr_H^G}$.

To see this suppose that $\ell(f) \geq p^t + 1$.  Then $f = \del^{p^t}(F)$ for some $F \in \field[V_n]$.
Therefore $f(e_i) =(\del^{p^t} F)(e_i) = ((\sigma-1)^{p^t}F)(e_i)=(\sigma^{p^t}(F) - F)(e_i)
= F(\sigma^{-p^t}(e_i)) - F(e_i)$.
Thus $f(e_i)=0$ if $e_i$ is fixed by $\sigma^{p^t}$, i.e., if $i \geq n-p^t+1$.
Therefore if $\ell(f) \geq p^t+1$ then $f$ vanishes on the set
$\mathrm{span}_\field\{e_{n-p^{t}+1},e_{n-p^{t}+2},\dots,e_{n-1},e_n\}$.
Hence if $\ell(f) \geq p^t+1$ then $f \in  \sqrt{\Im \Tr_H^G}$
\end{proof}

\begin{proposition}\label{basic decomp}
Suppose that $f$ is a non-zero homogeneous element of $\field[V_n]^G$. Then $\ell(f N^G(x_n)) = \ell(f)$.
\end{proposition}
\begin{proof}
Denote $N^G(x_n)$ by $N$.
Define $t$ such that $p^{t-1}<n\leq p^t$ (the case $n=1$ is trivial). Then the leading term of $N$ is $x_n^{p^t}$.
Let $\field[V_n]^\flat$ denote the span of the monomials in $\field[V_n]$ which, as polynomials in $x_n$, have degree less than $p^t$.
The fact that $x_n\not\in\del(\field[V_n])$ means that $\field[V_n]^\flat$ is a $\field G$-submodule of $\field[V_n]$.
For an arbitrary polynomial $h\in\field[V_n]$, viewing $h$ as a polynomial in $x_n$ and dividing by $N^G(x_n)$ gives
$h=qN+r$ for unique $r\in \field[V_n]^\flat$ and $q\in\field[V_n]$. This gives the $\field G$-module decomposition
$\field[V_n]=N\field[V_n]\oplus \field[V_n]^\flat$ (compare with \cite[Lemma~2.9]{hughes-kemper} and \cite[\S~2]{Shank-Wehlau:noether}).  As noted above $\ell(Nf) \geq \ell(f)$.  Suppose
$Nf = \del^t(F)$ and write $F=NF_1 + F_0$ with $F_0 \in \field[V_n]^\flat$.  Then
$N f = \del^t(N F_1 + F_0) = N \del^t(F_1) + \del^t(F_0) = N \del^t(F_1)$ and thus
$f = \del^t(F_1)$.  This shows that $\ell(f) \geq \ell(Nf)$.
\end{proof}

\section{Computing $\fvpp^\zps$}\label{comp_inv_sec}

In this section we use the ladder technique described in \cite[\S 7]{Shank-Wehlau:cmipg} to prove Theorem~\ref{invthm}.
We use the notation
$$G:=\langle\sigma\rangle\cong \zps,\ L:=\langle\sigma^p\rangle\cong \zp,\ {\rm and}\ Q:=G/L=\langle\Qg\rangle\cong \zp.$$
Note that $\deg\left(N^G\left(x_p\right)\right)=p$, $\deg\left(N^G\left(x_{p+1}\right)\right)=p^2$ and $Tr_L^G(x_{p})=x_1$.

The action of $L$ on $V^*_{p+1}$ is given by $\sigma^p(x_{p+1})=x_{p+1}+x_1$ and
$\sigma^p(x_i)=x_i$ for $i\leq p$. Thus as $L$-modules, $\fvpp\cong \field[V_2\oplus (p-1)V_1]$.
Therefore
$\fvpp^L\cong\field[\Ln,x_p,\ldots,x_1]$ with $\Ln=x^p_{p+1}-x_1^{p-1}x_{p+1}$. The action of $Q$ on $\fvpp^L$ is given by
$\Qg(\Ln)=\Ln +x^p_p-x_1^{p-1}x_p$ and $\Qg(x_i)=\sigma(x_i)$ for $i=1,2,\dots,p$.
Define
$$A:=\field[z_p,\ldots,z_1,X_p,\ldots,X_1]$$
with $\deg(z_i)=p$ and $\deg(X_i)=1$.
Further define an algebra homomorphism
$\pi:A\to\fvpp^L$ by $\pi(z_i)=x_{i+1}^p-x_1^{p-1}x_{i+1}$ and $\pi(X_i)=x_i$.
Note that $\pi$ is a degree preserving surjection with $\pi(z_p)=\Ln$. Further note that the
kernel of $\pi$ is the ideal
$$I:=\left(z_{p-1}-\left(X_p^p-X_1^{p-1}X_p\right),\ldots,z_1-\left(X_2^p-X_1^{p-1}X_2\right)\right)A.$$
Define an action,  by algebra automorphisms, of $Q$ on $A$ by taking $\Qg(z_i)=z_i+z_{i-1}$ and $\Qg(X_i)=X_i+X_{i-1}$ for $i>1$,
$\Qg(z_1)=z_1$ and $\Qg(X_1)=X_1$. Thus as $\field Q$-modules $A\cong \field[2V_p]$ and $\pi$ is a map of $\field Q$-modules.

The short exact sequence of $\field Q$-modules, $0\to I \to A \lato{\pi} \fvpp^L\to 0$, gives a long exact sequence on group cohomology
$$0\to I^Q\to A^Q \to \left(\fvpp^L\right)^Q\to \qhom 1 I \to \qhom 1 A \to \cdots .$$
We will show that the inclusion of $I$ into $A$ induces an injection of $\qhom 1 I$ into $\qhom 1 A$. Thus $\pi$ restricts to a surjection from
$A^Q$ to $\left(\fvpp^L\right)^Q=\fvpp^G$. Since $2V_p$ is a permutation representation of $Q$, after a suitable change of basis, $Q$ acts on $A$
by permuting the variables. Using the permutation basis, the orbit sums of monomials form a vector space basis
for $A^Q$. Since $Q\cong \zp$, these orbits are of size $p$ or size $1$. The orbits of size $p$ span
projective $\field Q$-module summands of $\fvpp^L$ while the orbits of size $1$ span trivial summands.
It easy to see that, in the original basis,
the orbits of size $1$ are polynomials in
$N^Q(z_p)$ and $N^Q(X_p)$, while the orbit sums coming from orbits of size $p$ are elements in the
image of the transfer. Thus $A^Q$ is generated by $N^Q(z_p)$, $N^Q(X_p)$ and elements from $\tr^Q(A)$, giving Theorem~\ref{invthm}.


The rest of this section is devoted to completing the proof of Theorem~\ref{invthm} by showing
that  the inclusion of $I$ into $A$ induces an injection of $\qhom 1 I $ into $\qhom 1 A $
(see Theorem~\ref{inclusion}~(b)).  We start by describing $H^*(Q,A)$.

\begin{proposition}\label{cohprop}
(a) $H^2(Q,A)=A^Q/\tr^Q(A)\cong \field[ N^Q(X_p), N^Q(z_p)]$.\\
(b) $\qhom 1 A $ is a principal $A^Q$-module with annihilator given by $\tr^Q(A)$.
\end{proposition}

\begin{proof} It follows from the discussion above
that, as an $\field Q$-module, $A$ consists of projective
summands and trivial summands  with the trivial summands
spanned by the monomials in $N^Q(X_p)$ and $N^Q(z_p)$. The projective summands do not contribute to the cohomology.
The trivial summands contribute non-zero classes to both the first and second cohomology.
\end{proof}

Note that, although $\qhom 1 A $ does not have
a multiplicative structure, it is isomorphic to  $\field[ N^Q(X_p), N^Q(z_p)]$ as an $A^Q$-module.

To compute $H^*(Q,I)$, we start by resolving $I$ as an $A-\field Q$-module using a
Koszul resolution (see, for example, \cite[\S 1.6]{bruns-herzog}).
Observe that $I$ is generated by an $A$-regular sequence
of length $p-1$. Furthermore, these generators span the degree $p$ homogeneous component, $I_p$, of $I$ and,
as a $\field Q$-module, $I_p\cong V_{p-1}$.
Let $\mu$ denote the $Q$-equivariant map from $V_{p-1}\otimes A$ to $A$ given by identifying elements of $V_{p-1}$ with elements of $I_p$
and then using the
multiplication in $A$.
Let $\Lambda^i(V_{p-1})$ denote
the $i^{\rm th}$ exterior power of $V_{p-1}$.
Define $\zeta^i:\Lambda^i(V_{p-1})\to \Lambda^{i-1}(V_{p-1})\otimes V_{p-1}$ by
$$\zeta^i(v_1\wedge v_2\wedge\cdots\wedge v_i)=\sum_{j=1}^i(-1)^{i-j}\left(v_1\wedge \cdots\wedge\widehat{v_j}\wedge \cdots\wedge v_i\right)\otimes v_j$$
for all $v_1,v_2,\dots,v_i\in V_{p-1}$.
 Define $F^{-i}:=\Lambda^i(V_{p-1})\otimes A$ for $i=1,2,\dots,p-1$ and
define $\rho^{-i}: F^{-i}\to F^{-i+1}$ to be $(1_{\Lambda^{i-1}(V_{p-1})}\otimes \mu) \circ(\zeta^i\otimes 1_A)$. This gives the following
sequence of $A-\field Q$-modules:
$$0\to F^{1-p} \lato{\rho^{1-p}} F^{2-p}\lato{\rho^{2-p}} \cdots \lato{\rho^{-3}} F^{-2}\lato{\rho^{-2}} F^{-1}\lato{\mu} I \to 0.$$
Since the generators of $I$ form a regular $A$-sequence, it follows from \cite[Corollary~1.6.14]{bruns-herzog} that this sequence is exact.
For $i>0$, define $K^{-i}$ to be the kernel of the map $\rho^{-i}:F^{-i}\to F^{-i+1}$. For convenience,
we define $K^0:=I$, $K^1:=A/I$, $K^a:=0$ for $a>1$, $F^1:=A/I$ and $F^a:=0$ for $a>1$. Using the exactness of the resolution,
we get a series of short exact sequences $0\to K^{-i}\to F^{-i}\to K^{-i+1}\to 0$. For each of these short exact sequences,
we apply group cohomology to get
a long exact sequence:
\begin{eqnarray*}
0\to \qhom 0 {K^{-i}} \to \qhom 0 {F^{-i}} &\to& \qhom 0 {K^{-i+1}} \to \qhom 1 {K^{-i}}\\ & \to & \qhom 1 {F^{-i}} \to \qhom 1 {K^{-i+1}} \to \ldots .
\end{eqnarray*}
Defining $D^{a,b}:=\qhom b {K^a}$ and $E^{a,b}:=\qhom b {F^a}$ gives a bigraded exact couple which leads to a spectral sequence.
This is essentially the construction given at the end of \cite[Ch XI, \S 5]{maclane}.
We will use this spectral sequence to describe $\qhom 1 I$. The following series of lemmas lead to a description of $\qhom b {F^a}$.

\begin{lemma}\label{splitlem} The map $\zeta^{i+1}$ is an isomorphism of $\field Q$-modules from $\Lambda^{i+1}(V_{p-1})$ to a direct summand of
$\Lambda^{i}(V_{p-1})\otimes V_{p-1}$.
\end{lemma}
\begin{proof}
It is clear that $\zeta^{i+1}$ followed by the natural projection from  $\Lambda^{i}(V_{p-1})\otimes V_{p-1}$ to
$\Lambda^{i+1}(V_{p-1})$ is $i+1$ times the identity map on $\Lambda^{i+1}(V_{p-1})$. Since $i+1\leq p-1$, this
is an isomorphism.
\end{proof}

\begin{lemma}\label{prodlem} $V_{p-1}\otimes V_{p-1}\cong V_1\oplus (p-2)V_p$.
\end{lemma}
\begin{proof} If the representation ring $R_{\field \zp}$ is extended by adjoining an element $\am$
satisfying $V_2=\am+\am^{-1}$, then by \cite[Lemma 2.3]{hughes-kemper},
$V_n=\frac{\am^n-\am^{-n}}{\am-\am^{-1}}$.
Thus, in the augmented representation ring,
\begin{eqnarray*}
V_{p-1}^2-V_p\cdot V_{p-2}&=&
\frac{(\am^{p-1}-\am^{-(p-1)})^2-(\am^p-\am^{-p})(\am^{p-2}-\am^{-(p-2)}))}{(\am-\am^{-1})^2}\\
&=&\frac{\am^{2p-2}-2+\am^{-(2p-2)}-(\am^{2p-2}-\am^2-\am^{-2}+\am^{-(2p-2)})}{\am^2-2+\am^{-2}}\\
&=&1=V_1.
\end{eqnarray*}
Therefore, $V_{p-1}\otimes V_{p-1}\cong V_p\otimes V_{p-2} \oplus V_1$.
It follows from \cite[Ch.~II \S 7 Lemma~4]{alperin} that $V_p\otimes V_i\cong i V_p$.
Thus $V_{p-1}\otimes V_{p-1}\cong (p-2)V_p\oplus V_1$, as required.
\end{proof}

\begin{lemma}\label{struclem} For $i$ even,
$$\Lambda^i(V_{p-1})\cong V_1\oplus \frac{1}{p}\left({p-1 \choose i}-1\right)V_p$$
and for $i$ odd,
$$ \Lambda^i(V_{p-1})\cong V_{p-1}\oplus \frac{1}{p}\left({p-1 \choose i}-(p-1)\right)V_p.$$
\end{lemma}
\begin{proof} First observe that $\dim_\field(\Lambda^i(V_{p-1})={p-1\choose i}\equiv (-1)^i \pmod p$.
Therefore the dimensions are correct. Thus it follows from Lemma~\ref{splitlem}
that $\Lambda^i(V_{p-1})$ is a non-projective summand of $\Lambda^{i-1}(V_{p-1})\otimes V_{p-1}$.
Also note that the result is true for $i=0$ and $i=1$.
We proceed by induction. Suppose the result holds for $i$. For $i$ even this gives,
\begin{eqnarray*}
\Lambda^i(V_{p-1})\otimes V_{p-1}& \cong&
        V_1\otimes V_{p-1}\oplus  \frac{1}{p}\left({p-1 \choose i}-1)\right)V_p\otimes V_{p-1}\\
      & \cong & V_{p-1} \oplus \frac{p-1}{p}\left({p-1 \choose i}-1)\right)V_p
\end{eqnarray*}
and therefore, $\Lambda^i(V_{p-1})\otimes V_{p-1}$ is isomorphic to $V_{p-1}$ plus projective modules.
Hence $\Lambda^{i+1}(V_{p-1})$ is isomorphic to $V_{p-1}$ plus projective modules.
For $i$ odd,
\begin{eqnarray*}
\Lambda^i(V_{p-1})\otimes V_{p-1} &\cong&
   V_{p-1}\otimes V_{p-1}\oplus  \frac{1}{p}\left({p-1 \choose i}-(p-1)\right)V_p\otimes V_{p-1}\\
   &\cong&V_{p-1}\otimes V_{p-1}\oplus  \frac{p-1}{p}\left({p-1 \choose i}-(p-1)\right)V_p
\end{eqnarray*}
and, therefore, $\Lambda^i(V_{p-1})\otimes V_{p-1}$ is isomorphic to $V_{p-1}\otimes V_{p-1}$ plus projective modules.
From Lemma~\ref{prodlem}, $V_{p-1}\otimes V_{p-1}$ is isomorphic to $V_1$ plus projective modules.
Hence $\Lambda^{i+1}(V_{p-1})$ is isomorphic to $V_1$ plus projective modules.
\end{proof}

\begin{lemma}\label{princlem} For $a\leq 0$ and $b>0$, $\qhom b {F^a}$ is a principal $A^Q$-module with annihilator $\tr^Q(A)$.
\end{lemma}
\begin{proof} As an $\field Q$-module, $A$ is a direct sum of projective summands with socles contained in $\tr^Q(A)$
and one dimensional summands spanned by monomials in $N^Q(X_p)$ and $N^Q(z_p)$. It follows from Lemma~\ref{struclem}
that $\Lambda^{-a}(V_{p-1})$ contains a a single non-projective summand. Note that projective summands do not contribute to the cohomology.
Further note that for any module $M$ and projective module $P$, $M\otimes P$ is projective.
Thus $\qhom b {\Lambda^{-a}(V_{p-1}})$ is a one dimensional vector space and
 $\qhom b {F^a}\cong \qhom b {\Lambda^{-a}(V_{p-1}}\otimes A) \cong \qhom b A$.
The result follows from Proposition~\ref{cohprop}
\end{proof}

The following lemma is a preliminary step in evaluating $d^{a,b}:E^{a,b}\to E^{a+1,b}$
for $a<0$ and $b>0$.

\begin{lemma}\label{Iplem} The inclusion of $I_p$ into $A_p$ induces an injection from
$\qhom 1 {I_p}$ to $\qhom 1 {A_p}$ and the zero map from $\qhom 2 {I_p}$ to $\qhom 2 {A_p}$.
\end{lemma}
\begin{proof} To see that the inclusion induces an injection from  $\qhom 1 {I_p}$ to $\qhom 1 {A_p}$, first
note that $\del$ is a twisted derivation and that $\del(f)$, for $f$ a generator of $A$, lies in
${\rm Span}_\field(z_{p-1},\ldots,z_1,X_{p-1}, \ldots, X_1)$.
Therefore $\del(A)$ is contained in the ideal $(z_{p-1},\ldots,z_1,X_{p-1},\ldots, X_1)A$.
As an $\field Q$-module, $I_p$ is isomorphic to $V_{p-1}$ with generator
$r:=z_{p-1}-\left(X_p^p-X_1^{p-1}X_p\right)$.
Thus $\qhom 1 {I_p}$ is a one dimensional vector space with $r$
representing a non-zero cohomology class. Since  $r$
does not lie in the ideal
$(z_{p-1},\ldots,z_1,X_{p-1},\ldots, X_1)A$,
this element does not lie in $\del(A)$ and, therefore, represents a non-zero class in $\qhom 1 {A_p}$.

To see that inclusion induces the zero  map from $\qhom 2 {I_p}$ to $\qhom 2 {A_p}$, observe that
$I_p\cong V_{p-1}$ and $A_p$ is isomorphic to $V_1$ plus projectives. Thus the map on cohomology is determined by
a $\zp$-equivariant map from $V_{p-1}$ to $V_1$, and all such
maps induce the zero map in second cohomology.
\end{proof}

For $a<0$, the first differential in the spectral sequence is the map on cohomology $d^{a,b}:\qhom b {F^a}\to \qhom b {F^{a+1}}$ induced by
$\rho^{a}:F^a \to F^{a+1}$

\begin{theorem}\label{diffthm}
For  $b>0$ and $a<0$, $d^{a,b}:\qhom b {F^a}\to \qhom b {F^{a+1}}$ is an isomorphism if $a$ and $b$ have the same parity and
zero if $a$ and $b$ have different parities.
\end{theorem}

\begin{proof} From Lemma~\ref{princlem}, $\qhom b {F^a}$ is a principal $A^Q$-module. Since $d^{a,b}$ is an $A^Q$-module map,
it is sufficient to evaluate $d^{a,b}$ on a generator. Thus, using the definition of $\rho^{a}$, we see that $d^{a,b}$ is
determined by the composition
$$\Lambda^{-a}(V_{p-1})\lato{\zeta^{-a}}\Lambda^{-a-1}(V_{p-1})\otimes V_{p-1}\lato{\cong}\Lambda^{-a-1}(V_{p-1})\otimes I_p\lato{\subset}
\Lambda^{-a-1}(V_{p-1})\otimes A.$$
It follows from Lemmas~\ref{splitlem} and \ref{struclem} that $\zeta^{-a}$ induces an isomorphism in cohomology.
Thus $d^{a,b}$ is determined by the inclusion of $\Lambda^{-a-1}(V_{p-1})\otimes I_p$ into $\Lambda^{-a-1}(V_{p-1})\otimes A.$

For $a$ odd, using Lemma~\ref{struclem}, $\Lambda^{-a-1}(V_{p-1})$ is isomorphic to $V_1$ plus projectives. Therefore, in this case,
$d^{a,b}$ is induced by the inclusion of $I_p$ into $A_p$. Thus, using Lemma~\ref{Iplem}, if $a$  and $b$ are both odd then $d^{a,b}$ is
injective and if $a$ is odd and $b$ is even  then $d^{a,b}$ is zero.

For $a$ even, using Lemma~\ref{struclem}, $\Lambda^{-a-1}(V_{p-1})$ is isomorphic to $V_{p-1}$ plus projectives.
Therefore, in this case, $d^{a,b}$ is induced by the inclusion of $V_{p-1}\otimes I_p$ into
$V_{p-1}\otimes A_p \cong V_{p-1}\otimes V_{p-1}$.
By Lemma~\ref{prodlem}, $V_{p-1}\otimes V_{p-1}$ is isomorphic to $V_1$ plus projectives. Since $A_p$ is isomorphic to
$V_1$ plus projectives, $V_{p-1}\otimes A_p$ is isomorphic to $V_{p-1}$ plus projectives.
Thus $d^{a,b}$ is determined by the composition
$$V_1\to V_{p-1}\otimes I_{p} \to V_{p-1}\otimes A_p \to V_{p-1}.$$
This map clearly induces the zero map from $\qhom 1 {V_1}$ to $\qhom 1 {V_{p-1}}$. Thus for $a$ even and $b$ odd, $d^{a,b}=0$.
To show that $d^{a,b}$ is injective for $a$ even and $b$ even, we need to show that the given map from $V_1$ to $V_{p-1}$ is non-zero.
It follows from Lemma~\ref{Iplem}, that for the purposes of computing cohomology,
the inclusion of $I_p$ into $A_p$ is
the injection of $V_{p-1}$ into $V_1\oplus V_p$ taking $e'$ to $(e'',\del(e))$
where $e, e'$ and $e''$ denote elements which generate the cyclic $G$-modules
$V_p$, $V_{p-1}$ and $V_1$ respectively.
The cokernel of this map is isomorphic to
$V_2$. Tensoring over $\field$
is exact so we have a short exact sequence
$$0\to V_{p-1}\otimes V_{p-1} \to V_{p-1}\otimes (V_1\oplus V_p)\to V_{p-1}\otimes V_2 \to 0.$$
This gives rise to a long exact sequence in cohomology. Recall that $V_{p-1}\otimes V_2\cong V_{p-2}\oplus V_p$
(see, for example, \cite[Ch.II~\S~7~Lemma~5]{alperin}).
Thus, modulo projectives, the sequence is $V_1 \to V_{p-1} \to V_{p-2}$. This can only give a long exact sequence on cohomology
if the the map from $V_1$  to $V_{p-1}$ is non-zero.
\end{proof}

\begin{corollary}\label{E2cor} For $b>0$ and $a<0$, the spectral sequence satisfies
$$E_2^{a,b}=\left\{
\begin{array}{ll}
\field & {\rm if}\ a=1-p\ {\rm and}\ b\ {\rm odd};\\
 0 & {\rm otherwise}.
\end{array}
\right.
$$
\end{corollary}

It follows from Theorem~\ref{diffthm} that $\rho^{-1}$ induces an isomorphism
from $\qhom 1 {F^{-1}}$ to $\qhom 1 {A}$. This map factors through $\qhom 1 I$ with
the first map in the factorisation induced by $\mu$ and the second induced by inclusion.
Thus to complete the proof of Theorem~\ref{invthm}, it is sufficient to show the following.

\begin{lemma} The map
$\mu$ induces an epimorphism from $\qhom 1 {F^{-1}}$ to $\qhom 1 I$.
\end{lemma}

\begin{proof}
Denote by $\partial^{a,b}$ the connecting homomorphism from $\qhom b {K^a}$ to $\qhom {b+1} {K^{a-1}}$
and define a filtration on $\qhom 1 I = \qhom 1 {K^0}$ by
$${\mathcal F}_t:={\rm kernel}( \partial^{-t,t+1}\circ \partial^{-t+1,t}\circ \cdots \circ \partial^{-2,3} \circ \partial^{-1,2}  \circ \partial^{0,1}).$$
Since $\partial^{1-p,p}=0$, we have ${\mathcal F}_{p-1}=\qhom 1 I$.
Using the long exact sequence in cohomology coming from $0\to K^{-1}\to F^{-1}\to K^0\to 0$, we see that ${\mathcal F}_0$ is the image of
$\qhom 1 {F^{-1}}$ in $\qhom 1 I$. We will prove the lemma by showing ${\mathcal F}_0={\mathcal F}_{p-1}$.

Using the definition a derived couple (see, for example, \cite[Ch. XI, \S 5]{maclane}), we have
$D_{t+2}^{a-t-1,b+t+1}=\partial^{a-t,b+t}\circ \cdots  \circ \partial^{a-1,b+1}  \circ \partial^{a,b}(D^{a,b})$.
If $x\in {\mathcal F}_t \setminus {\mathcal F}_{t-1}$ then
$\partial^{-t+1,t}\circ\cdots \circ \partial^{0,1}(x)$ is a non-zero element of $D_{t+1}^{-t,t+1}$ which lifts to a  non-zero element of
$E_{t+1}^{-(t+1),t+1}$. However, it follows from Corollary~\ref{E2cor} that $E_2^{-(t+1),t+1}=0$ for $t\geq 1$.
Thus $E_{t+1}^{-(t+1),t+1}=0$ for $t\geq 1$. Therefore ${\mathcal F}_t={\mathcal F}_0$ for all $t \geq 1$.
\end{proof}

These calculations give the following.

\begin{theorem}\label{inclusion}
 (a) $\qhom 1 I $ is a principal $A^Q$-module with generator represented by $z_{p-1}-\left(X_p^p-X_1^{p-1}X_p\right)$
and  annihilator $\tr^Q(A)$.\\
(b) The inclusion of $I$ into $A$ induces an $A^Q$-module monomorphism of
$\qhom 1 I $ to $\qhom 1 A $ taking $[z_{p-1}-\left(X_p^p-X_1^{p-1}X_p\right)]$
to $-[N^Q(X_p)]$.
\end{theorem}

This completes the proof of Theorem~\ref{invthm}.

\section{The Noether number of $V_{p+1}$}\label{noeth_sec}

In this section we use the description of $\fvpp^\zps$
given in Theorem~\ref{invthm} to prove the following.

\begin{theorem} For $p>2$, the Noether number of $V_{p+1}$ is $p^2+p-3$.
\end{theorem}

\begin{remark} A Magma \cite{magma} calculation shows that for $p=2$, the Noether number of $V_3$ is $p^2=4$.
\end{remark}

For the remainder of this section we will assume that $p \geq 3$.
We continue to use the notation described at the beginning of
Section~\ref{comp_inv_sec}. Define $M:=N^G(x_p)$ and $N:=N^G(x_{p+1})$.
The theorem is an immediate consequence of the following two lemmas.

\begin{lemma}  
    The Noether number of $V_{p+1}$ is less than or equal to $p^2+p-3$. 
\end{lemma}

\begin{proof}
Let ${\mathcal H}$ denote the ideal in $\fvpp^L$ generated by the homogeneous $G$-invariants of positive degree,
i.e., ${\mathcal H}=\fvpp^G_+\cdot\fvpp^L$.
Thus $\fvpp^L/{\mathcal H}$ is a finite dimensional graded algebra, the ring of relative coinvariants.
Let ${\mathcal B}$ denote the set of elements of $\fvpp^L$ of the form $\gamma\cdot x_p^j\cdot N^L(x_{p+1})^k$,
with  $\gamma$ a monomial in $\{x_1,\ldots,x_{p-1}\}$ of degree at most $p-2$ and
$j,k<p$. The methods of Section~3 of \cite{fssw} show that ${\mathcal B}$ projects to a spanning set
in  $\fvpp^L/{\mathcal H}$. Therefore $(p-1)p+(p-1)+p-2=p^2+p-3$ is an upper bound on
the top degree of the relative coinvariants and $\Tr_L^G({\mathcal B})$ is a generating set for the ideal
$\Im \Tr_L^G$. By Theorem~\ref{invthm}, $\fvpp^G$ is generated by $N$, $M$ and elements from $\Im \Tr_L^G$.
Thus $p^2+p-3$ is an upper bound for the Noether number.
\end{proof}

\begin{lemma} 
The polynomial $\Tr_L^G\left(\left(N^L\left(x_{p+1}\right)x_p\right)^{p-1}x_{p-1}^{p-2}\right)$
is indecomposable in $\fvpp^G$.  In particular the Noether Number of $V_{p+1}$ is at least $p^2+p-3$.
\end{lemma}

\begin{proof}
Define $w:=N^L(x_{p+1})$ and $z:=\Tr_L^G\left(w^{p-1}x_p^{p-1}x_{p-1}^{p-2}\right)$.
Suppose, by way of contradiction, that $z=f_1h_1+\cdots +f_sh_s$
where $f_i$ and $h_i$ are homogeneous positive degree elements of $\fvpp^G$. The degree of $z$ as a polynomial in
$x_{p+1}$ is less than $p^2$. Thus $N$ does not appear in the decomposition.

We use the graded reverse lexicographic term order with $x_1<x_2<\cdots <x_{p+1}$ and denote the leading monomial
of an element $f\in\fvpp$ by $\lm(f)$. It is easy to see that $\lm(M)=x_p^p$.
An elementary calculation gives $\lm(z)=x_p^{p^2-1}x_{p-1}^{p-2}$.
By relabelling if necessary, we may assume $\lm(f_ih_i)\geq \lm(f_{i+1}h_{i+1})$. Thus, either $\lm(f_1h_1)=\lm(z)$ or
$\lm(f_1h_1)=\lm(f_2h_2)>\lm(z)$.
Without loss of generality, we may assume $f_1h_1=cM^m\alpha$, where $c\in\field$ and $\alpha$ is a 
(non-constant) product of elements from $\Tr_L^G({\mathcal B})$.

Let $\pi$ denote the projection $$\pi:\fvpp\to\fvpp/(x_1,\ldots,x_{p-2},x_{p-1}^{p-1})\fvpp.$$
For convenience, write $f\equiv h$ if $\pi(f)=\pi(h)$.
Observe that $\pi(z)\not=0$, $\pi(w)\equiv x_{p+1}^p$ and $\pi(M)\equiv x_p^p$.
Furthermore, the restriction of $\pi$ to $\fvpp^L$ commutes with the action of $Q=G/L$.
Thus $\pi\Tr_L^G({\mathcal B})=\Tr^Q\pi({\mathcal B})$. If
$\beta\in \Tr_L^G({\mathcal B})$ with $\pi(\beta)\not=0$, then $\beta=\Tr_L^G(w^kx_p^jx_{p-1}^{\ell})$
and $\pi(\beta)\equiv x_{p-1}^{\ell}\Tr^Q(w^kx_p^j)$. Summing over the action of $Q$ gives
\begin{eqnarray*}
\Tr^Q(w^kx_p^j)&\equiv& \sum_{\lambda\in\fp}(x_{p+1}^p+\lambda x_p^p)^k(x_p+\lambda x_{p-1})^j\\
              &\equiv& \sum_{\lambda\in\fp}\left(\sum_{t=0}^k\binomial{k}{t}\lambda^t x_{p+1}^{p(k-t)}x_p^{tp}\right)
                       \left(\sum_{r=0}^j\binomial{j}{r}\lambda^r x_p^{j-r}x_{p-1}^r\right)\\
              &\equiv& \sum_{r=0}^j\sum_{t=0}^k \left(\sum_{\lambda\in\fp}\lambda^{r+t}\right)\binomial{k}{t} \binomial{j}{r}
                        x_{p+1}^{p(k-t)}x_p^{tp+j-r}x_{p-1}^r.
\end{eqnarray*}
Recall that $\sum_{\lambda\in\fp}\lambda^i=0$ unless $i$ is a non-zero multiple of $p-1$, in which case the sum is $-1$.
Therefore $\Tr^Q(w^kx_p^j)\equiv 0$ if $j+k<p-1$.  Moreover,
if $p-1\leq j+k<2p-2$, we take $t+r=p-1$ to get
\begin{eqnarray*}
\Tr^Q(w^kx_p^j) &\equiv& -\sum_{t=1}^k \binomial{k}{t}\binomial {j}{p-1-t} x_{p+1}^{p(k-t)}x_p^{pt+j+t-(p-1)}x_{p-1}^{p-1-t}\\
               &\equiv& -\sum_{r=p-1-k}^j \binomial{k}{p-1-r}\binomial {j}{r} x_{p+1}^{p(k+r-(p-1))}x_p^{p(p-1-r)+j-r}x_{p-1}^r\\
               &\equiv& -\binomial{j}{p-1-k} x_p^{k(p+1)+j-(p-1)}x_{p-1}^{p-1-k}+x_{p-1}^{p-k}F
\end{eqnarray*}
with $F\in\field[x_{p+1},x_p,x_{p-1}]$. Since $j\leq p-1$ and $k\leq p-1$, we have $j+k\geq 2p-2$ only when $j=p-1$
and $k=p-1$. In this case, there is one additional term, $-x_p^{p^2-p}x_{p-1}^{p-1}\equiv 0$. Since the monomials of
degree $kp+j$ taken to zero by $\pi$ are less than $x_p^{kp+j-p-2}x_{p-1}^{p-2}$, we have, for $k>0$,
$$\lm\left(\Tr_L^G\left(w^kx_p^j\right)\right)=x_p^{k(p+1)+j-(p-1)}x_{p-1}^{p-1-k}.$$

Assume, by way of contradiction, that $\alpha$ is the product of at least two factors,
 say $\alpha=\beta_1\beta_2\cdots\beta_d$ with  $\beta_i\in\Tr_L^G({\mathcal B})$.
Since  $\lm(cM^m\alpha)\geq\lm(z)=x_p^{p^2-1}x_{p-1}^{p-2}$, we have $\lm(\alpha)\geq x_p^{p^2-mp-1}x_{p-1}^{p-2}$.
Therefore, since we are using the graded reverse lexicographic order, $\pi\lm(\alpha)\not=0$.
Furthermore, $\lm(\beta_i)$ divides $\lm(\alpha)$. Thus $\pi(\beta_i)\not=0$ giving
$\beta_i=\Tr_L^G(w^{k_i}x_p^{j_i}x_{p-1}^{\ell_i})$ with $j_i+k_i\geq p-1$. Using the formulae above
gives $$\lm(\beta_1\beta_2)=x_p^{(p+1)(k_1+k_2)+j_1+j_2-2(p-1)}x_{p-1}^{2(p-1)-k_1-k_2+\ell_1+\ell_2}.$$
Again, using  $\lm(cM^m\alpha)\geq\lm(z)$ gives $2(p-1)-k_1-k_2+\ell_1+\ell_2\leq p-2$ which simplifies to
$k_1+k_2\geq p+\ell_1+\ell_2\geq p$. However,
${\rm deg}(\beta_1\beta_2)=p(k_1+k_2)+j_1+j_2+\ell_1+\ell_2\leq p^2+p-3$, giving $k_1+k_2\leq p$. Therefore $k_1+k_2=p$.
Furthermore, adding the inequalities $j_i+k_i\geq p-1$ gives $j_1+j_2+k_1+k_2=j_1+j_2+p\geq 2(p-1)$ which simplifies to
$j_1+j_2\geq p-2$. Thus ${\rm deg}(\beta_1\beta_2)\geq p(k_1+k_2)+j_1+j_2\geq p^2+p-2>  {\rm deg}(cM^m\alpha)$, giving a contradiction.
Thus we must have that $\alpha$ is an element of $\Tr_L^G({\mathcal B})$.

It remains to consider the case $f_1h_1=cM^m\alpha$ with $\alpha\in \Tr_L^G({\mathcal B})$ and $m>0$.
As above, $\pi\lm(\alpha)\not=0$ gives $\alpha=\Tr_L^G(w^kx_p^jx_{p-1}^{\ell})$ with
$k+j\geq p-1$ and $\ell+p-1+k\leq p-2$. The degree constraint gives $p(m+k)+j+\ell=p^2+p-3$.
Since $\alpha\in\Tr_L^G({\mathcal B})$, we have $j,k\leq p-1$ and $\ell \leq p-2$.
Thus $j+\ell\leq 2p-3$. Therefore, either
$m+k=p$ and $j+\ell=p-3$ or $m+k=p-1$ and $j+\ell=2p-3$.

We first consider the case $m+k=p-1$ and $j+\ell=2p-3$. Since $j\leq p-1$ and $\ell\leq p-2$, we have $j=p-1$ and $\ell=p-2$.
Using the above formula for $\pi\Tr^Q(w^kx_p^j)$, with $m>0$,gives
$$\pi\left(M^m\Tr^G_L\left(w^{p-1-m}x_p^{p-1}x_{p-2}^{p-2}\right)\right)
\equiv \binomial{p-1}{m} x_p^{p^2-1-m} x_{p-1}^{p-2+m}\equiv 0.$$
Thus $\lm\left(M^m\Tr^G_L\left(w^{p-1-m}x_p^{p-1}x_{p-2}^{p-2}\right)\right)<\lm(z)$.

This leaves the case
$m+k=p$ and $j+\ell=p-3$.
Again using the formula for $\pi\Tr^Q(w^kx_p^j)$ gives
\begin{eqnarray*}
\lm\left(M^m\Tr^G_L\left(w^{p-m}x_p^jx_{p-1}^{p-3-j}\right)\right)&=&x_p^{mp+kp+j-(p-1-k)}x_{p-1}^{p-1-k+\ell}\\
&=&x_p^{p^2+j-m+1}x_{p-1}^{p-4+m-j}.
\end{eqnarray*}
However, $j+k\geq p-1$ gives $m-j\leq 1$. Therefore $\lm(cM^m\alpha)\not=\lm(z)$. However, it is possible to choose
$\alpha$ and $m$ so that $\lm(cM^m\alpha)=x_p^{p^2+p-3-s}x_{p-1}^s$ for $s=0,\ldots,p-3$. Note that $s=p-4+m-j=\ell+m-1$.
Since $m\geq 1$, $s=0$ occurs only when $m=1$ and $\ell=0$. In general, we may take $m=1,2,\ldots,s+1$ and $\ell=s+1-m$.
Define $T_{m,s}:=M^m\Tr_L^G(w^{p-m}x_p^{p-4+m-s}x_{p-1}^{s+1-m})$. To complete the proof of the lemma, it is sufficient to show
that no linear combination of elements of ${\mathcal S}:=\{T_{m,s}\mid s=0,\ldots,p-3,\ m=1,2,\ldots,s+1\}$
has lead monomial $\lm(z)=x_p^{p^2-1}x_{p-1}^{p-2}$. Our argument is essentially
Gauss-Jordan elimination applied to $\pi{\mathcal S}$.

Using the above formula for
$\pi\Tr^Q(w^kx_p^j)$ gives
$$T_{m,s}\equiv -\sum_{r=m-1}^{p-4+m-s}\binomial{p-m}{p-1-r}\binomial{p-4+m-s}{r}x_{p+1}^{p(r-m+1)}x_{p-1}^{r+s+1-m}x_p^*.$$
Reindexing with $i=r+s+1-m$ gives
\begin{eqnarray*}
T_{m,s}&\equiv&-\sum_{i=s}^{p-2} \binomial{p-m}{p-i+s-m}\binomial{p-4+m-s}{i-s+m-1}x_{p+1}^{p(i-s)}x_{p-1}^ix_p^*\\
       &\equiv&-\sum_{i=s}^{p-2} \binomial{p-m}{i-s}\binomial{p-4+m-s}{i-s+m-1}x_{p+1}^{p(i-s)}x_{p-1}^ix_p^*.
\end{eqnarray*}
Note that in a field of characteristic $p$,
$\binomial{p-a}{b}=(-1)^a\binomial{a-1+b}{a-1}$. Thus
$$T_{m,s}\equiv (-1)^{s+1}\sum_{i=s}^{p-2}\binomial{m-1+i-s}{m-1}\binomial{i+2}{s-m+3}x_{p+1}^{p(i-s)}x_{p-1}^ix_p^*.$$
A simple calculation confirms $\binomial{a}{c}\binomial{a+b}{b}=\binomial{b+c}{c}\binomial{a+b}{b+c}$, giving
\begin{eqnarray*}
T_{m,s}&\equiv& (-1)^{s+1}\sum_{i=s}^{p-2}\binomial{s+2}{m-1}\binomial{i+2}{s+2}x_{p+1}^{p(i-s)}x_{p-1}^ix_p^{p^2+p-3-i-p(i-s)}\\
       &\equiv& (-1)^{s+1}\binomial{s+2}{m-1}\sum_{i=s}^{p-2}\binomial{i+2}{s+2}x_{p+1}^{p(i-s)}x_{p-1}^ix_p^{p^2+p-3-i-p(i-s)}.
\end{eqnarray*}
Therefore
$$(-1)^{s+1}T_{m,s}\binomial{s+2}{m-1}^{-1}\equiv\sum_{i=s}^{p-2}\binomial{i+2}{s+2}x_{p+1}^{p(i-s)}x_{p-1}^ix_p^{p^2+p-3-i-p(i-s)}$$
is independent of $m$. Thus $\{\pi(T_{1,s})\mid s=0,\ldots,p-3\}$ is a basis for $\pi{\mathcal S}$. Since
$\{\lm(T_{1,s})\mid s=0,\ldots,p-3\}=\{x_p^{p^2+p-3-s}x_{p-1}^s\mid s=0,\ldots, p-3\}$, no linear combination of elements of
${\mathcal S}$ has lead monomial $x_p^{p^2-1}x_{p-1}^{p-2}$.

  The final assertion of the lemma follows from the fact that $\deg z = p^2+p-3$.
\end{proof}

\section{Decomposing $\fvpp$}\label{decomp sec}

The main goal of this section is to describe the $\field G$-module decomposition of $\fvpp$.
We do this by considering a basis for $\fvpp^G_d$ which is compatible with the length
filtration.  By Theorem~\ref{invthm}, $\fvpp^G$ is generated by $M:=N^G(x_p)$ and $N:=N^G(x_{p+1})$
together with elements of the relative transfer, $\Tr_L^G(\fvpp^L)$. To identify the summands
occurring in the decomposition, we need to determine the lengths of the basis elements.

Suppose $f\in \fvpp^G$.
It follows from Lemma~\ref{length lemma} that
if $\ell(f) \geq 2$ then $f \in \sqrt{\Im \Tr_L^G} = ((x_1,x_2,\dots,x_p)\fvpp)\cap\fvpp^G$, and
if $\ell(f) \geq p+1$ then $f \in \sqrt{\Im \Tr^G} = ((x_1)\fvpp)\cap\fvpp^G$.
Furthermore, since $\Tr_L^G=\del^{p-1}$, if $\ell(f)\geq p$ then $f\in\Im\Tr_L^G$.
Since $N = N^G(x_{p+1})$ has lead term\footnote{Use any monomial order with $x_1<x_2<\cdots <x_{p+1}$.}
$x_{p+1}^{p^2}$, we have that
$N \notin (x_1,x_2,\dots,x_p)\fvpp$ and thus $\ell(N) = 1$.
Similarly since the lead term of $M = N^G(x_p)$ is $x_p^p$, we have that
$M \notin (x_1)\fvpp$ and thus $\ell(M) \leq p$.

\begin{lemma}\label{delqp_lem}  Let $1 \leq q < p$.   Then
 $$\del^{qp}(x_{p+1}^i) =
   \begin{cases}
     0,              & \text{if } i < q\\
     q!x_1^q,  & \text{if } i=q\\
     x_1^q h \text{\quad for some } h \in \fvpp,& \text{if } i \geq {q+1}.
  \end{cases}$$
In particular, $\del^{qp}(x_{p+1}^i) \in (x_1^q)\fvpp$ for all $i \geq 0$.
\end{lemma}
\begin{proof}
  We consider $\del^{qp}(x_{p+1}^i)$ using induction on $q$.
For $q=1$ we have
\begin{eqnarray*}
\del^p(x_{p+1}^i) &=& (x_{p+1}+x_1)^i - x_{p+1}^i
= \sum_{j=0}^{i-1} \binomial{i}{j}x_1^{i-j}x_{p+1}^j\\
&=&\begin{cases}
     0,               & \text{if } i =0\\
     x_1,          & \text{if } i=1\\
     x_1  \sum_{j=0}^{i-1} \binomial{i}{j}x_1^{i-j-1}x_{p+1}^j,    & \text{if } i \geq 2.
  \end{cases}
\end{eqnarray*}
Now take $q+1 \geq 2$.  Then
\begin{eqnarray*}
\del^{(q+1)p}(x_{p+1}^i) &=& \del^{qp}(\del^p(x_{p+1}^i))\\
 &=& \del^{qp}(\sum_{j=0}^{i-1} \binomial{i}{j}x_1^{i-j}x_{p+1}^j)\\
&=& \sum_{j=0}^{i-1} \binomial{i}{j}x_1^{i-j}\del^{qp}(x_{p+1}^j)\ .
\end{eqnarray*}
By induction this gives
\begin{eqnarray*}
\del^{(q+1)p}(x_{p+1}^i) & = &
\sum_{j=q}^{i-1} \binomial{i}{j}x_1^{i-j}\del^{qp}(x_{p+1}^j)\\
&=& \begin{cases}
        0,                                                                   & \text{if } i-1 < q\\
        \binomial{q+1}{q}x_1\del^{qp}(x_{p+1}^q),   & \text{if } i-1=q\\
        x_1 \sum_{j=q}^{i-1} \binomial{i}{j}x_1^{i-j-1}\del^{qp}(x_{p+1}^j), & \text{if } i-1 > q
\end{cases}\\
&=& \begin{cases}
        0,                                 & \text{if } i < q+1\\
        (q+1) x_1 q! x_1^q,   & \text{if } i=q+1\\
         x_1 \sum_{j=q}^{i-1} x_1^q h_j\text{\quad where } h_j \in \fvpp,  &  \text{if } i \geq {q+1}.
\end{cases}\\
&= &\begin{cases}
        0,                                 & \text{if } i < q+1\\
        (q+1)!  x_1^{q+1},   & \text{if } i=q+1\\
         x_1^{q+1} h \text{\quad for some } h \in \fvpp, &  \text{if } i \geq {q+1}.
\end{cases}
\end{eqnarray*}
\end{proof}

\begin{proposition}\label{del^p prop}
  Let $f$ be a non-zero element of $\fvpp^G$ and let $1 \leq q < p$.
Then $\ell(f) \geq qp+1 \iff f \in \Im \del^{qp} \iff f \in (x_1^q)\fvpp^G$.
\end{proposition}
\begin{proof}
  The first equivalence is just the definition of length.  For the second equivalence,
first suppose that $ f = \del^{qp}(F) \in \Im \del^{qp}$.
Write $F = \sum_{i=0}^r f_i x_{p+1}^i$ where each $f_i \in \field[x_1,x_2,\dots,x_p]$.
Then $f = \del^{qp}(F) = \sum_{i=0}^r f_i \del^{qp}(x_{p+1}^i) \in (x_1^q)\fvpp^G$
by the previous lemma.
  Conversely, suppose that $f \in (x_1^q)\fvpp^G$ and write $f = x_1^q f'$ where
$f' \in \fvpp^G$.   Then $$\del^{qp}\left(\frac{x_{p+1}^q f'}{q!}\right) = \frac{f'}{q!} \del^{qp}(x_{p+1}^q)= \frac{f'}{q!}q! x_1^q = f.$$
\end{proof}

\begin{proposition} \label{stretch prop}
  Let $f$ be a non-zero element of $\fvpp^G$ and write $f=x_1^q f'$ where $x_1$ does not divide $f'$.
If $q \geq p$ then $\ell(f) = p^2$.  Otherwise $\ell(f) = qp + \ell(f')$.
\end{proposition}
\begin{proof}
Applying Lemma~\ref{delqp_lem} gives $\del^{(p-1)p}(x_{p+1}^{p-1})=(p-1)!x_1^{p-1}=-x^p_1$.
Furthermore, $\del^p(x_p)=0$.
Thus
\begin{eqnarray*}
\del^{p^2-1}(x_{p+1}^{p-1}x_p)&=&\del^{p-1}\left((\del^p)^{p-1}(x_{p+1}^{p-1}x_p)\right)
=\del^{p-1}(x_p\del^{(p-1)p}(x_{p+1}^{p-1}))\\
&=&-\del^{p-1}(x_px_1^{p-1})=-x_1^{p-1}\del^{p-1}(x_p)=-x_1^p.
\end{eqnarray*}
Therefore $\del^{p^2-1}(-x_{p+1}^{p-1}x_pf')=x_1^pf'$. This implies that if $q \geq p$ then $\ell(f) = p^2$.
Suppose then that $q < p$.  By Proposition~\ref{del^p prop}, we have
$qp \leq \ell(f)-1$.  Since $x_1^{q+1}$ does not divide $f$,  Proposition~\ref{del^p prop}
also implies that $\ell(f)-1 < (q+1)p$.
Write $\ell(f)-1 = qp + r$ where $0 \leq r \leq p-1$ and define $s := \ell(f')-1$.  Since
$x_1$ does not divide $f'$, Lemma~\ref{length lemma} implies that $0 \leq s \leq p-1$.
We will show that $r=s$.

  Clearly there exists $F\in \fvpp$ such that $f = \del^{qp+r}(F)$.
Therefore $f=\del^r(\del^{qp}(F)) = \del^r(x_1^q F')$ for some $F' \in \fvpp$.
Hence $x_1^q f' = f = x_1^q \del^r(F')$ and therefore $f' = \del^r(F')$.
Hence $s+1 = \ell(f') \geq r+1$.

  Conversely we may write $f' = \del^s(F'')$ for some $F'' \in \fvpp$.
 Since $s \leq p-1$ we have $\del^p(F'') = \del^{p-s}(\del^s(F'')) = \del^{p-s}(f') = 0$.
This shows that $F'' \in \fvpp^L$.
Thus $\del^{qp+s}(x_{p+1}^q F'') = (\del^s (\del^p)^q)(x_{p+1}^q F'') = \del^s(q! x_1^q F'')
= q! x_1^q \del^s(F'') = q! x_1^q f' = q! f$ where $q! \neq 0$ since $q < p$.
This shows that $f \in \Im \del^{qp+s}$ and thus $qp+r+1 = \ell(f) \geq qp+s +1$.
Therefore $r=s$ as required.
\end{proof}

\begin{proposition}\label{nondiv_tr_prop}
  Let $f$ be a non-zero element in the image of the relative transfer, $\Tr_L^G(\fvpp^L)$.
Suppose that $x_1$ does not divide $f$.   Then $\ell(f) = p$.
\end{proposition}
\begin{proof}
  Since $x_1$ does not divide $f$,  Lemma~\ref{length lemma} implies that
$\ell(f) \leq p$.  Conversely $\Tr_L^G = 1 + \sigma + \sigma^2 + \dots + \sigma^{p-1} = \del^{p-1}$.
Thus the hypothesis that $f \in \Im \Tr_L^G$ implies that $\ell(f) \geq p$.
\end{proof}

\begin{remark} Since elements in $\tr^G(\field[V_{p+1}])$ have length $p^2$, it is clear that
$\field[V_{p+1}]^G$ is generated by $N$, $M$ and elements from  $\Im \Tr_L^G\setminus \Im \tr^G$.
\end{remark}

\begin{proposition}\label{M length}
  $\ell(M^j) = j+1$ for all $j=0,1,\dots,p-1$.  In particular $M^j$ lies in the image of the relative transfer,
$\Tr_L^G(\fvpp^L)$, if and only if $j \geq p-1$.
\end{proposition}

\begin{proof}
From Theorem~\ref{invthm}, $\fvpp^G$ is generated by $M$, $N$ and elements from $\Im \Tr_L^G$.
Note that $\deg(M)=p$ and $\deg(N)=p^2$. Thus for $d<p^2$, if $p$ does not divide $d$, we have $\fvpp^G_d=\Tr_L^G(\fvpp^L_d)$
and, if $d=ip$ with $i<p$, $\fvpp^G_d=\field\cdot M^i+\Tr_L^G(\fvpp^L_{ip})$.
Fix $j\in\{1,2,\ldots,p-1\}$.
Choose a basis, $\B$, for $\fvpp_{jp}^G$ so that $\B$ is compatible with the length filtration.
Applying Proposition~\ref{decomp prop} gives a decomposition
$$\fvpp_{jp}=\bigoplus_{\alpha\in \B} V(\alpha)\ .$$

Suppose $f\in \F{p}^G(\fvpp)\cap\B$. Then $f\in \Tr_L^G(\fvpp^L)$. If $x_1$ does not divide $f$,
then by Proposition~\ref{nondiv_tr_prop}, $\ell(f)=p$. Suppose $x_1$ does divide $f$. Write $f=x_1^qf'$
where $x_1$ does not divide $f'$. If $f'\in \Tr_L^G(\fvpp^L)$, then by Proposition~\ref{stretch prop}
and Proposition~\ref{nondiv_tr_prop}, $\ell(f)$ is a multiple of $p$. If $f'\not\in \Tr_L^G(\fvpp^L)$,
then $f'=cM^i+f''$ for some non-zero $c\in\field$ and some $f''\in \Tr_L^G(\fvpp^L)$. Thus $\deg(f')=ip$
and $q=(j-i)p$. Thus $q>p$ and by Proposition~\ref{stretch prop}, $\ell(f)$ is $p^2$.
Hence, for any $f\in \F{p}^G(\fvpp)\cap\B$, $\ell(f)$ is a multiple of $p$.
Therefore $p$ divides the dimension of $$\bigoplus_{\alpha\in \F{p}^G(\fvpp)\cap \B} V(\alpha).$$

Since $\dim \fvpp_{pj} = \binomial{p+pj}{pj}$, Lucas'~Lemma (see, for example, \cite{fine})
implies that $\dim \fvpp_{pj} \equiv \binom{j+1}{j} \binom{0}{0} \pmod p$.
Thus $\dim \fvpp_{pj} \equiv j+1 \pmod p$.   This shows that $M^j \notin \Tr_L^G(\fvpp^L)$
for all $j \leq p-2$ and that $M^{p-1} \in  \Tr_L^G(\fvpp^L)$.  Furthermore, for
$1 \leq j < p-1$ we have $\ell(M^j) < p$ and $\ell(M^j) \equiv j+1 \pmod p$ and thus
$\ell(M^j) = j+1$.  Since $\fvpp_0 \cong \field$, it is clear that $\ell(M^0) = 1$.
\end{proof}

Define $\fvpp^\flat$ to be the span of the monomials in $\fvpp$ which,
as polynomials in $x_{p+1}$, have degree less than $p^2$.
It follows from the proof of Proposition~\ref{basic decomp} that, as $\field G$-modules,
$\fvpp=N\fvpp\oplus \fvpp^\flat$. Thus a decomposition of $\fvpp^\flat$ gives a decomposition of
$\fvpp$.  Therefore the following theorem implies Theorem~\ref{decompthm0}.

\begin{theorem}\label{decompthm}
(i) For $d<p^2-p$, an $\field G$-module decomposition of
$\fvpp^\flat_d$ includes precisely one non-induced indecomposable summand.
Divide $p$ into $d$ to get $d=bp+c$ with $0\leq c<p$. The non-induced
summand is isomorphic to $V_{cp+b+1}$ and the decomposition may be chosen so that
the socle of the non-induced indecomposable is spanned by $x_1^cM^b$.\\
(ii) For $d\geq p^2-p$, $\fvpp^\flat_d$ is a direct sum of indecomposable induced modules.
\end{theorem}

\begin{proof}
Fix $d$ and choose a basis, $\B$, for $(\fvpp^{\flat}_d)^G$, so that $\B$ is compatible with the length filtration.
Applying Proposition~\ref{decomp prop} gives a decomposition
$$\fvpp^{\flat}_d=\bigoplus_{\alpha\in \B} V(\alpha)$$
with $V(\alpha)\cong V_{\ell(\alpha)}$.
Write $\alpha=x_1^i\alpha'$ where $x_1$ does not divide $\alpha'$.
If $i\geq p$, then by Proposition~\ref{stretch prop}, $\ell(\alpha)=p^2$
and $V(\alpha)$ is projective. Suppose $i<p$. If $\alpha'\in\Im\Tr_L^G$, then
by Proposition~\ref{nondiv_tr_prop}, $\ell(\alpha')=p$. Thus, using Proposition~\ref{stretch prop},
$\ell(\alpha)=ip+p$ and $V(\alpha)$ is an induced module. Suppose $\alpha'\not\in\Im\Tr_L^G$.
Then $\alpha'=kM^j+h$ where $j<p-1$, $k$ is a non-zero element of $\field$ and $h\in\Im\Tr_L^G$.
It follows from Proposition~\ref{M length} that $\ell(\alpha')=j+1$. Applying Proposition~\ref{stretch prop}
gives $\ell(\alpha)= pi+j+1$. This last case is the only way in which a non-induced summand can appear in the decomposition.
Note that in this case, $d=pj+i$ with $0\leq i<p$ and $j<p-1$, giving $d\leq (p-2)p+(p-1)=p^2-p-1$, $i=c$ and $j=b$.
Suppose, by way of contradiction, that $\alpha_1,\alpha_2\in\B$ are distinct elements both having length $cp+b+1$. Then  $\alpha_1=x_1^c(k_1M^b+h_1)$ and $\alpha_2=x_1^c(k_2M^b+h_2)$
with $k_i\in \field\setminus\{0\}$ and $h_i\in\Im\Tr_L^G$. From Proposition~\ref{stretch prop},
 $\ell(k_2\alpha_1-k_1\alpha_2)\geq cp+p>cp+b+1$ contradicting the fact that $\B$ is compatible with the length filtration.
\end{proof}

\begin{remark} The strong form of the Hughes-Kemper Periodicity Conjecture \cite[Conjecture~4.6]{hughes-kemper}
states that for $p^{m-1}<n\leq p^m$ and $d>p^m-n$, $\field[V_n]^{\flat}_d$ is induced. The preceding Theorem verifies
the conjecture for $n=p+1$.
\end{remark}

\section{Appendix}
  In this appendix we will derive generating functions which give the multiplicities of the indecomposable 
  $\field G$-modules as summands in $\field[V_{p+1}]_n$.
    We also derive the Hilbert series for $\field[V_{p+1}]^G$.  

  Throughout the appendix we will write $n = \alpha p^2 + \beta p +\gamma$ where
  $0 \leq \beta,  \gamma \leq p-1$.
If $\beta \neq p-1$ then by Theorem~\ref{decompthm0}
we know that $\field[V_{p+1}]_n$ contains exactly one non-induced
summand, $V_{d(n)}$ where $d(n) = \gamma p + \beta + 1$. For
convenience we will also define $d(n) = \gamma p + \beta + 1=
(\gamma+1)p$ when $\beta = p-1$.

  Define integer valued functions $a_1(n), a_2(n), \dots, a_p(n)$ by
    \begin{equation}\label{decomp}
    \field[V_{p+1}]_n \cong V_{d(n)} \oplus a_1(n)\,V_p \oplus
    a_2(n)\,V_{2p} \oplus \dots \oplus a_p(n)\, V_{p^2}\ .
    \end{equation}
  By Propositions~\ref{del^p prop} and \ref{stretch prop}, an invariant $f$ spans the socle of a copy of $V_{ip}$ where $p>i \geq 2$,
  if and only if $f=x_1h$ where the invariant $h$ spans the socle of a copy of $V_{(i-1)p}$.  Clearly if $n=\deg(f)$ then
  $\deg(h)=n-1\geq 0$.  This means that for all $2 \leq i \leq p-1$ we have
\begin{equation}\label{recursive1}
  a_i(n) =
   \begin{cases}
       0,               & \text{if } n=0,1;\\
       a_{i-1}(n-1),    & \text{if } n \geq 1.
   \end{cases}
\end{equation}
Similarly, Proposition~\ref{stretch prop} with $q =p-1$ and $q=p$ combined with
Theorem ~\ref{decompthm0} implies that
 \begin{equation}\label{recursive2}
 a_p(n) =
   \begin{cases}
       0,                            & \text{if } n=0,1;\\
       a_p(n-1) + a_{p-1}(n-1),      & \text{if $p$ does not divide $n$};\\
       a_p(n-1) + a_{p-1}(n-1) + 1,& \text{if $p$ divides $n$ and } n \neq 0.
   \end{cases}
 \end{equation}
   Furthermore comparing dimensions in the
   decomposition~(\ref{decomp})
yields the equation:
   \begin{equation}\label{dim eqn}
      \binomial{n+p}{p} = d(n) + p a_1(n) + 2p a_2(n) + \dots + p^2 a_p(n)\ .
   \end{equation}

Introduce the generating functions:
\begin{eqnarray*}
  D(x) &=& \sum_{n=0}^\infty d(n) x^n \\
  A_i(x) &=& \sum_{n=0}^\infty a_i(n) x^n \text{ for } i=1,2,\dots,p.
\end{eqnarray*}
  In terms of these generating functions, the above recursive
  conditions, (\ref{recursive1}) and (\ref{recursive2}), become:
  \begin{eqnarray*}
     A_i(x) & = & \sum_{n=0}^\infty a_i(n) x^n \\
            & = & a_i(0) + x\sum_{n=1}^\infty a_{i-1}(n-1)x^{n-1}\\
            & = & x A_{i-1}\qquad  (\text{for } i=2,3,\dots,p-1)
  \end{eqnarray*}
  and
 \begin{eqnarray}\label{Recursive2}
     A_p(x) & = & \sum_{n=0}^\infty a_p(n) x^n \nonumber\\ 
            & = & a_p(0) + x\sum_{n=1}^\infty \left(a_p\left(n-1\right) + a_{p-1}\left(n-1\right) + \delta_n^0\right)
             x^{n-1} \nonumber\\
             & &\qquad \text{ where }\delta^0_n = \begin{cases}
                                               1,& \text{if } n \equiv 0 \pmod p;\\
                                               0,& \text{otherwise},
                                        \end{cases}\nonumber\\
            & = & x A_p(x) + x A_{p-1}(x) + \sum_{n=1}^\infty x^{np}\nonumber\\
            & = & x A_p(x) + x A_{p-1}(x) + \frac{x^p}{1-x^p}\ .
   \end{eqnarray}
Again using the generating functions, the dimension equation~(\ref{dim eqn}) becomes
    $$\frac{1}{(1-x)^{p+1}} = D(x) + p A_1(x) + 2p A_2(x) + \dots + p^2 A_p(x) \ .$$
  Substituting $A_2(x) = x A_1(x)$, $A_3(x) = x^2 A_1(x)$,
  $\ldots$, $A_{p-1}(x) = x^{p-2} A_1(x)$, we are left with the following two
  equations in $A_1$ and $A_p$:
\begin{eqnarray*}\label{system1}
   A_p(x) &=& x A_p(x) + x^{p-1} A_1(x) + \frac{x^p}{1-x^p}\\
   p A_1(x) + 2p x A_1(x)  
    + &\dots& + (p^2-p) x^{p-2} A_1(x) + p^2 A_p(x) +D(x) = \frac{1}{(1-x)^{p+1}}.\nonumber
\end{eqnarray*}
Collecting terms this system 
 becomes:
\begin{eqnarray*}\label{system2}
   -x^{p-1} A_1(x) + (1-x) A_p(x) &=& \frac{x^p}{1-x^p}\\
   p(1 + 2x + 3x^2 + \dots + (p-1)x^{p-2}) A_1(x) + p^2 A_p(x) &=&  - D(x) + \frac{1}{(1-x)^{p+1}}.\nonumber
\end{eqnarray*}

  Note that, as is easily verified by integration, we have
  \begin{equation}\label{powereqn}
  1 + 2x + 3x^2 + \dots + mx^{m-1} = \frac{1-(m+1) x^m + mx^{m+1}}{(1-x)^2}.
  \end{equation}
  Thus the above system of equations 
   becomes
\begin{eqnarray}\label{system3}
   -x^{p-1} A_1(x) + (1-x) A_p(x) &=& \frac{x^p}{1-x^p}\\
   p\left(\frac{1-px^{p-1}+(p-1)x^p}{(1-x)^2}\right) A_1(x) + p^2 A_p(x) &=&  - D(x) + \frac{1}{(1-x)^{p+1}}.\nonumber
\end{eqnarray}


\smallskip\noindent
  Solving for $A_1$ and $A_p$ yields
\begin{eqnarray*}\label{solution}
\lefteqn
{A_1(x) = \left(\frac{-p^2x^p}{1-x^p} + \frac{1}{(1-x)^p}-(1-x)D(x)\right) \left(\frac{1-x}{p(1-x^p)}\right)}\\
 A_p(x) &=& \biggl(\frac{1-p x^{p-1}+(p-1)x^p}{(1-x)^2}\left(\frac{px^p}{1-x^p}\right)+
                  \frac{x^{p-1}}{(1-x)^{p+1}} - x^{p-1}D(x)\biggr)   \left(\frac{1-x}{p(1-x^p)}\right)\ .
\nonumber
\end{eqnarray*}

Thus a closed form for $D(x)$ will yield closed forms for
$A_1(x)$ and $A_p(x)$.  To obtain a closed expression for $D(x)$ we observe that the sequence $\{d(n)\}_{n=0}^\infty$ is the sum of two periodic sequences, one of period $p$ and one of period $p^2$.
From this using Equation~(\ref{powereqn}) twice we get
\begin{eqnarray*}
D(x)&=& \left(\sum_{\gamma=0}^{p-1}p\gamma x^\gamma\right)\frac{1}{1-x^p}
              + \left(\sum_{\beta=0}^{p-1} \sum_{\gamma=0}^{p-1}(\beta+1)x^{p\beta+\gamma}\right)\frac{1}{1-x^{p^2}} \\
 &=& px\left(\sum_{\gamma=0}^{p-1}\gamma x^{\gamma-1}\right)\frac{1}{1-x^p}               + \left(\sum_{\beta=0}^{p-1} (\beta+1)x^{p\beta}\right)
                 \left(\sum_{\gamma=0}^{p-1}x^\gamma\right)\frac{1}{1-x^{p^2}}\\
 &=& px\left(\frac{1-px^{p-1}+(p-1)x^p}{(1-x)^2}\right)\frac{1}{1-x^p}\\
& &\quad        + \left(\frac{1-(p+1)(x^p)^p+p(x^p)^{p+1}}{(1-x^p)^2}\right)
        \left(\frac{1-x^p}{1-x}\right)\frac{1}{1-x^{p^2}}.
\end{eqnarray*}

%
%
\noindent  Substituting this expression into the expression for $A_p$ given above 
  and simplifying yields the following:
  \begin{eqnarray*}
  A_p(x) = \frac{1}{p(1-x^p)} 
       \left(\frac{x^{p-1}}{(1-x)^p} - \frac{x^{p-1}-(p+1)x^{p^2+p-1}+px^{p^2+2p-1}}{(1-x^p)(1-x^{p^2})}\right)
       \ .
  \end{eqnarray*}
Using this expression for $A_p(x)$ in (\ref{system3})
 gives
 \begin{eqnarray*}
  A_1(x) = \frac{-x}{1-x^p} + \frac{1-x}{p(1-x^p)} 
       \left(\frac{1}{(1-x)^p} - \frac{1-(p+1)x^{p^2}+px^{p^2+p}}{(1-x^p)(1-x^{p^2})}\right)
       \ .
  \end{eqnarray*}

\medskip
  In the above description, the summand $V_{d(n)}$ is sometimes an
  induced summand.  More precisely, this happens exactly when
  $\beta=p-1$.  Thus if we decompose the induced component
  \begin{equation*}
  (\field[V_{p+1}]_n)_{\text{induced}} \cong b_1(n)  V_p \oplus b_2(n) V_{2p} \oplus \dots \oplus b_p(n) V_{p^2}
  \end{equation*}
  we have
  $$b_i(n) =
  \begin{cases}
    a_i(n) + 1, & \text{if $\gamma=i-1$ and $\beta=p-1$};\\
    a_i(n),     &\text{ otherwise}.
  \end{cases}$$
Thus the generating function $B_i(x) = \sum_{n=0}^\infty b_i(n)
x^n$ is given by
$$B_i(x) = A_i(x) + \frac{x^{p^2-p+i-1}}{1-x^{p^2}} = x^{i-1}A_1(x)+ \frac{x^{p^2-p+i-1}}{1-x^{p^2}}.$$

  Note that the Hilbert series of the ring of invariants
  $\field[V_{p+1}]^{{\bf Z}/p^2}$ is given by
$\mathcal{H}(\field[V_{p+1}]^{{\bf Z}/p^2},x) = \frac{1}{1-x} + \sum_{i=1}^p A_i(x)$.
Repeatedly using the recursive equation for $A_p(n)$ (\ref{Recursive2}) we obtain
\begin{eqnarray*}
A_p & = & xA_p + xA_{p-1} + \frac{x^p}{1-x^p}\\
        & = & x( xA_p +  xA_{p-1} + \frac{x^p}{1-x^p}) +  xA_{p-1} + \frac{x^p}{1-x^p}\\
        & = & x^2A_p +  (x^2 + x)A_{p-1} + (x+1)\frac{x^p}{1-x^p}\\
        & = & x^2(xA_p +  xA_{p-1} + \frac{x^p}{1-x^p}) +  (x^2 + x)A_{p-1} + (x+1)\frac{x^p}{1-x^p}\\
         & = & x^3A_p +  (x^3+x^2 + x)A_{p-1} + (x^2+x+1)\frac{x^p}{1-x^p}\\
         & &\qquad\qquad\qquad\qquad\vdots\\
         & = & x^{p-1}A_p +  (x^{p-1}+x^{p-2} + \dots + x)A_{p-1} +  (x^{p-2}+x^{p-3} + \dots + 1)\frac{x^p}{1-x^p}\\
         & = & x^{p-1}\left(A_p +  (1+x^{-1} + \dots + x^{-(p-2)})A_{p-1}\right) + \left(\frac{1-x^{p-1}}{1-x}\right) \frac{x^p}{1-x^p}\\
    & = & x^{p-1}\left(A_p +  (A_{p-1}+A_{p-2} + \dots + A_1) + \left(\frac{1}{1-x}- \frac{1}{1-x^p}\right)\right)\\
    &=& x^{p-1}\left(\mathcal{H}(\field[V_{p+1}]^{{\bf Z}/p^2},x) - \frac{1}{1-x^p}\right).\\
\end{eqnarray*}

 \noindent Therefore
   \begin{eqnarray*}
   \mathcal{H}(\field[V_{p+1}]^{{\bf Z}/p^2},x) &=& \frac{1}{x^{p-1}}A_p(x) + \frac{1}{1-x^p}\\
   &=& \frac{1}{p(1-x^p)}
       \left(\frac{1}{(1-x)^p} + \frac{(p-1)-px^p+x^{p^2}}{(1-x^p)(1-x^{p^2})}\right)
   \ .
   \end{eqnarray*}

\ifx\undefined\bysame
\newcommand{\bysame}{\leavevmode\hbox to3em{\hrulefill}\,}
\fi

\end{document}